\documentclass[a4paper,12pt,twoside,leqno]{article}

\usepackage{amssymb} 
\usepackage{amsfonts}
\usepackage{amsmath}
\usepackage[utf8x]{inputenc}
\usepackage{tikz,relsize} 
\usepackage[cmtip,arrow]{xy}
\usepackage{pb-diagram,pb-xy}

\usepackage{amsthm}
\usepackage{hyperref}

\def\ns#1{\mathbb{#1}}
\def\N{\ns{N}}
\def\Z{\ns{Z}}
\def\Q{\ns{Q}}
\def\R{\ns{R}}

\DeclareMathOperator{\Isom}{Isom}
\DeclareMathOperator{\Aff}{Aff}
\DeclareMathOperator{\Aut}{Aut}
\DeclareMathOperator{\Inn}{Inn}
\DeclareMathOperator{\Out}{Out}
\DeclareMathOperator{\GL}{GL}

\numberwithin{equation}{section}

\theoremstyle{plain}
\newtheorem{thm}{Theorem}[section]
\newtheorem{lm}[thm]{Lemma}
\newtheorem{cor}[thm]{Corollary}
\newtheorem{prop}[thm]{Proposition}

\theoremstyle{definition}

\newtheorem{ex}{Example}[section]

\theoremstyle{remark}

\author{R. Lutowski, A. Szczepa\'nski
\thanks{Both authors are supported by the Polish National Science Center grant 2013/09/B/ST1/04125.}}

\title{Crystallographic groups with trivial center and outer automorphism group}

\date{}

\begin{document}
\maketitle

\renewcommand{\thefootnote}{}
\renewcommand{\thefootnote}{\arabic{footnote}}
\setcounter{footnote}{0}

\begin{abstract}
Let $\Gamma$ be a crystallographic group of dimension $n,$ i.e. a discrete, cocompact
subgroup of $\Isom(\R^n)$ = $O(n)\ltimes\R^n.$ For any $n\geq 2,$ we construct a crystallographic
group with a trivial center and a trivial outer automorphism group.
\end{abstract}

\section{Introduction}
Let $\Gamma$ be a discrete, cocompact subgroup of $O(n)\ltimes\R^n$ = $\Isom(\R^n)$ i.e.
a crystallographic group. 
If $\Gamma$ is a torsion free group, then $M = \R^n/\Gamma$ is a flat manifold (that is a compact Riemannian manifold without
boundary with the sectional curvature $K_x = 0$ for any $x\in M$). Moreover $\pi_{1}(M) = \Gamma.$
\vskip 1mm
\noindent
In 2003 R. Waldm\"uller found
a torsion free crystallographic group $\Gamma\subset O(141)\ltimes\R^{141}$ (a flat manifold $M = \R^{141}/\Gamma$) with 
the following properties:
\vskip 1mm
$(i)$ $Z(\Gamma) = \{e\},$
\vskip 1mm
$(ii)$ $\Out(\Gamma) = \{e\},$
\vskip 1mm
\noindent
where $Z(\Gamma)$ is the center of the group $\Gamma,$ and $\Out(\Gamma) = \Aut(\Gamma)/\Inn(\Gamma)$ denotes the group
of outer automorphisms of $\Gamma$ (see \cite[Appendix C]{S} and \cite{wald}). 
Equivalently, $(i)$  means that the abelianization of $\Gamma$ is finite (the first Betti number of $M$ is equal to zero). 
Moreover, if both conditions $(i)$ and $(ii)$ are satisfied, then the group
of affine diffeomorphisms $\Aff(M)$ of the manifold $M$ is trivial
(see \cite{charlap} and \cite{S}).
\vskip 1mm
\noindent
We do not know if there exist such flat manifolds
in dimensions less than $141.$ For example in dimensions up to six
such Bieberbach groups do not exist.  
In this paper we are interested in the existence of not necessarily torsion free crystallographic groups with the above properties.
We shall prove that
for any $n\geq 2$ there exists a crystallographic group of dimension $n$ which
satisfies conditions $(i)$ and $(ii).$

The main motivation for us is the article \cite{BL} of M. Belolipetsky and A. Lubotzky.
For any $n\geq 3$ they found an infinite family of hyperbolic compact manifolds of dimension $n$ 
with the following property: for every manifold $M$ from this family, $\Out(\pi_1(M))$ = $\{e\}.$
Since the center of the fundamental group of a compact hyperbolic manifold
is trivial, the above result gives us an infinite family of groups which satisfy conditions $(i)$ and $(ii).$
The construction of the above hyperbolic examples uses the properties of simple Lie groups of $\R$-rank one
and, in particular, follows from the existence of non arithmetic lattices.
\vskip 1mm
\noindent
In our construction the most important are Bieberbach theorems, and specific properties of crystallographic groups.
\vskip 1mm
\section{Crystallographic groups with trivial center and outer automorphism group}
In this part we shall prove our main result. Let $\Gamma$ be a torsion free crystallographic group. 
From Bieberbach's theorems (see \cite[Chapter 2]{S}) we have a short exact sequence of groups
$$0\to\Z^n\to\Gamma\stackrel{p}{\rightarrow} G\to 0,$$
where
$\Z^n$ is a maximal abelian subgroup of $\Gamma$ and
$G$ is a finite group. Moreover, let $h_{\Gamma}:G\to \GL(n,\Z)$
be the integral holonomy representation defined by the formula
$$\forall_{g\in G} h_{\Gamma}(g)(e) = \bar{g}e\bar{g}^{-1},$$
where $\bar{g}\in\Gamma, p(\bar{g}) = g$ and $e\in\Z^n.$
Let
$$N = N_{\GL(n,\Z)}(h_{\Gamma}(G)) = \{X\in \GL(n,\Z)\mid \forall_{f\in h_{\Gamma}(G)}\hskip 2mm XfX^{-1}\in h_{\Gamma}(G)\}.$$
In the case when $Z(\Gamma) = \{e\},$ we have the following commutative diagram (\cite[p. 65-69]{S}) with exact rows and columns:
$$\begin{diagram}
\node{}\node{0}\arrow{s}\node{0}\arrow{s}\node{0}\arrow{s}\\
\node{0}\arrow{e}\node{\Z^n}\arrow{s}\arrow{e}\node{\Gamma}\arrow{s}\arrow{e}\node{G}\arrow{s,r}{h_\Gamma}\arrow{e}\node{0}\\
\node{0}\arrow{e}\node{Z^1(G,\Z^n)}\arrow{s}\arrow{e}\node{\Aut(\Gamma)}\arrow{s}\arrow{e,t}{F}\node{N_\alpha}\arrow{s}\arrow{e}\node{0}\\
\node{0}\arrow{e}\node{H^1(G,\Z^n)}\arrow{s}\arrow{e}\node{\Out(\Gamma)}\arrow{s}\arrow{e}\node{N_\alpha/G}\arrow{s}\arrow{e}\node{0}\\
\node{}\node{0}\node{0}\node{0}
\end{diagram}$$
\vskip 1mm
\begin{center}
Diagram 1
\end{center}
where $Z^1(G,\Z^n)$ is the group of 1-cocycles.
Moreover
$$N_{\alpha} = \{n\in N\mid n\ast\alpha =\alpha\},$$ and $\alpha\in H^2(G,\Z^n)$ is the cohomology
class of the first row of the diagram.
The action $\ast:N\times H^2(G,\Z^n)\to H^2(G,\Z^n)$ is defined by the formula
$$n\ast [a] = [n\ast a],$$
where $n\in N, a\in Z^2(G,\Z^n),\hskip 2mm [a]$ is the cohomology class of $a$ and  
$$\forall_{g_1,g_2\in G} \;  n\ast a(g_1,g_2) = n a(n^{-1}g_{1}n,n^{-1}g_{2}n).$$
\vskip 1mm
\noindent
We have the following proposition.
\begin{prop}
$\Aut(\Gamma)$ is a crystallographic group if and only if $\Out(\Gamma)$ is a finite group.
\end{prop}
\noindent
{\bf Proof:} 
We start with an observation that $Z^1(G,\Z^n)$ is a free abelian group of rank $n$ which is a faithful $N_{\alpha}$ module. 
First, assume that $\Aut(\Gamma)$ is a crystallographic group with the maximal abelian subgroup $M.$ 
From \cite[Proposition I.4.1]{charlap}, $M$ is the unique normal maximal abelian subgroup of $\Aut(\Gamma).$
Hence, $M = Z^1(G,\Z^n),$ and $\Out(\Gamma)$ is a finite group. The reverse implication is obvious.
This finishes the proof of the proposition.
\vskip 1mm
\hskip 120mm $\Box$
\vskip 1mm
\noindent
Let us formulate our main result.
\begin{thm}\label{main} For every $n\geq 2$ there exists a crystallographic group $\Gamma$ of dimension $n$
with $Z(\Gamma) = \Out(\Gamma) = \{e\}.$
\end{thm}
\vskip 1mm
\noindent
{\bf Proof:}
\vskip 5mm
\noindent
We shall need the following lemma.
\begin{lm}
Let $G, H$ be finite groups and $H\subset G\subset \GL(n,\Z).$ If
the group $N_{\GL(n,\Z)}(H)$ is finite, then $N_{\GL(n,\Z)}(G)$ is finite. 
\end{lm}
\noindent
{\bf Proof of Lemma:} From the assumption, $\Aut(H)$ and $\Aut(G)$ are finite. Moreover,
we have monomorphisms:
$$N_{\GL(n,\Z)}(H)/C_{\GL(n,\Z)}(H)\stackrel{\bar{\phi}}{\rightarrow}\Aut(H)$$
and
$$N_{\GL(n,\Z)}(G)/C_{\GL(n,\Z)}(G)\stackrel{\bar{\phi}}{\rightarrow}\Aut(G),$$
where $\bar{\phi}$ is induced by $\phi(s)(g)=sgs^{-1},g\in G, s\in \GL(n,\Z).$
Since
$C_{\GL(n,\Z)}(G)\subset C_{\GL(n,\Z)}(H),$ our Lemma is proved.
\vskip 1mm
\hskip 120mm $\Box$
\begin{cor}
If $\mid\Out(\Gamma)\mid < \infty,$ then $\mid\Out(\Aut(\Gamma))\mid < \infty.$ 
\end{cor}
\vskip 1mm
\hskip 120mm $\Box$
\begin{lm}
Assume $Z(\Gamma) = \{e\},$ then
\begin{enumerate}
\item $H^1(G,\Z^n)\simeq (\Q^{n}/\Z^{n})^{G} = H^{0}(G,\Q^{n}/\Z^{n});$
\item $Z^1(G,\Z^n)\simeq \{m\in\Q^n\mid\forall_{g\in G}\hskip 2mm gm-m \in \Z^n\} = A^{0}(\Gamma)$ as $N_{\alpha}$ modules;
\item $A(\Gamma) = N_{\Aff(\R^n)}(\Gamma) = \{a\in \Aff(\R^n)\mid\forall_{\gamma\in\Gamma} a\gamma a^{-1}\in\Gamma\}\simeq \Aut(\Gamma).$
\end{enumerate}
\end{lm}
\vskip 1mm
\hskip 120mm $\Box$
\vskip 1mm
\noindent
We have the following modification of the Diagram 2. 
$$\begin{diagram}
\node{}\node{0}\arrow{s}\node{0}\arrow{s}\node{0}\arrow{s}\\
\node{0}\arrow{e}\node{\Z^n}\arrow{s}\arrow{e}\node{\Gamma}\arrow{s}\arrow{e}\node{G}\arrow{s,r}{h_\Gamma}\arrow{e}\node{0}\\
\node{0}\arrow{e}\node{A^{0}(\Gamma)}\arrow{s}\arrow{e}\node{A(\Gamma)}\arrow{s}\arrow{e,t}{F}\node{N_\alpha}\arrow{s}\arrow{e}\node{0}\\
\node{0}\arrow{e}\node{(\Q^{n}/\Z^{n})^{G}}\arrow{s}\arrow{e}\node{\Out(\Gamma)}\arrow{s}\arrow{e}\node{N_\alpha/G}\arrow{s}\arrow{e}\node{0}\\
\node{}\node{0}\node{0}\node{0}
\end{diagram}$$
\vskip 1mm
\begin{center}
Diagram 2
\end{center}

\vskip 1mm
\noindent
Let $\Gamma$ be a crystallographic group of rank $n$ with trivial center and holonomy group $G.$
Moreover, assume that the group $H^1(G,\Z^n) = \{e\},$ and the group  
$\Out(\Gamma)$ is finite. Inductively, put $\Gamma_{0} = \Gamma$
and $\Gamma_{i+1} =$A($\Gamma_{i}$), for $i\geq 0.$
\begin{lm}\label{seq}
$\exists N$ such that $\Gamma_{N+1} = \Gamma_{N}.$
\end{lm}
\noindent
{\bf Proof:} We start from observations that for $i > 0,$ $\Gamma_i$ is a crystallographic group, 
$Z(\Gamma_i) = \{e\}$ and $M_0 = M_i,$ where $M_i = A^{0}(\Gamma_{i-1})\subset\Gamma_i$
is the maximal abelian normal subgroup (a subgroup of translations). Let $G_i =\Gamma_i/M_i.$ From definition we can consider
$(G_i)$ as a nondecreasing sequence of finite subgroups of $\GL(n,\Z).$ 
From Bieberbach theorems \cite[Chapter 2]{S} and from Diagrams 1 and 2, there is only a finite number of possibilities for $G_i.$
Hence $\exists N\in\N$ such that $\forall_{i > N}\hskip 2mm G_i = G_{N}.$
This finishes the proof.
\vskip 1mm
\hskip 120mm $\Box$
\begin{ex}
Let $\Gamma_1 = G_1\ltimes\Z^2$ be the crystallographic group of dimension 2 with holonomy group $G_1 = D_{12},$
where 
$$D_{12} = \text{gen}\left\{ 
\left [
\begin{matrix}
0 & -1\cr
1 & -1
\end{matrix}\right ],
\left [
\begin{matrix}
-1 & 0\cr
0 & -1
\end{matrix}\right ],
\left [
\begin{matrix}
0 & 1\cr
1 & 0
\end{matrix}\right ]\right\}
$$ 
\noindent
is the dihedral group of order 12.
Moreover, let $\Gamma_2 = G_2\ltimes\Z^3$ be the crystallographic group of dimension 3, with holonomy group $G_2 = S_4\times\Z_2$
generated by
matrices
$$\left [
\begin{matrix}
0 & 1 & 0\cr
0 & -1 & -1\cr
1 & 1 & 0
\end{matrix}
\right ],
\quad
\left [
\begin{matrix}
0 & 0 & 1\cr
0 & -1 & -1\cr
-1 & 0 & 1
\end{matrix}
\right ].$$
Here $S_4$ denotes the symmetric group on four letters.
\vskip 3mm
\noindent
For $i = 1,2$ we have
$$N_{\GL(n_{i},\Z)}(G_{i}) = G_{i}$$ and
$$H^1(G_i,\Z^{n_i}) = 0,$$
where $n_i$ is the rank of $\Gamma_i.$
Hence $A(\Gamma_i) = \Gamma_i,$ and $\Out(\Gamma_i) = \{e\},$ for $i = 1,2.$
\end{ex}
\vskip 5mm
\noindent
Now we are ready to finish the proof of Theorem \ref{main}.
The cases $n = 2,3$ are done in the above example.
Assume $n\geq 4.$ Let $n = 2k+3i,$ where $i\in\{0,1\}.$
Put $\Gamma' =\Gamma_{1}^{k}\times\Gamma_{2}^{i}.$
Then $\Gamma'$ is centerless and by \cite[Theorem 3.4]{lut2}
the bottom exact sequence of the Diagram 2 looks as follows
$$0\to 0\to \Out(\Gamma')\to S_k\to 0.$$
Hence, $\Gamma'$ satisfies the assumption of Lemma \ref{seq} and the sequence
$\Gamma_{0} = \Gamma',\hskip 2mm \Gamma_{i+1} =$ A($\Gamma_{i}$) stabilizes, i.e.,
$\exists N$ such that $\forall_{i\geq N}\hskip 3mm \Gamma_{i}$ = $\Gamma_{N}.$
Moreover,  $\Out(\Gamma_{N}) = \{e\}$ and $Z(\Gamma_{N}) = \{e\}.$
\vskip 1mm
\hskip 120mm $\Box$

\vskip 3mm
\noindent
Institute of Mathematics\\
University of Gda\'nsk\\
ul. Wita Stwosza 57\\
80-952 Gda\'nsk\\
Poland\\
E-mail: \texttt{rlutowsk@mat.ug.edu.pl}, \texttt{matas@univ.gda.pl}

\end{document}